\documentclass[12pt,english]{article}
\usepackage{fullpage}
\usepackage{babel}
\usepackage{exscale,xspace}
\usepackage{amsfonts,amssymb}

\usepackage[
  pdfauthor={Bernd C. Kellner},
  pdfkeywords={Kurepa's left factorial, subfactorial, derangements},
  pdftitle={Some remarks on Kurepa's left factorial},
  colorlinks,bookmarksnumbered]{hyperref}

\newcommand{\binom}[2]{ {#1 \choose #2} }
\newcommand{\pdiv}{\mid}

\newcommand{\notdiv}{\nmid}

\newcommand{\mod}[1]{({\rm mod\ } #1)}

\newcommand{\Bell}{\mathcal{B}}

\def\ord{\mathop{\rm ord}\nolimits}

\newcommand{\ZZ}{\mathbb{Z}}
\newcommand{\NN}{\mathbb{N}}

\newcommand{\qed}{\parbox{0cm}{}\hspace*{\fill} $\Box$}
\newcommand{\refeqn}[1]{(\ref{#1})}

\renewcommand{\theequation}{\arabic{section}.\arabic{equation}}
\newcommand{\settheequation}[1]{\renewcommand{\theequation}{#1}}
\newcommand{\restoretheequation}{\renewcommand{\theequation}%
{\arabic{section}.\arabic{equation}}\addtocounter{equation}{-1}}

\newtheorem{prop}{Proposition}[section]
\newtheorem{theorem}[prop]{Theorem}
\newtheorem{corl}[prop]{Corollary}
\newtheorem{lemma}[prop]{Lemma}
\newtheorem{conj}[prop]{Conjecture}

\newtheorem{thdefin}[prop]{Definition}
\newtheorem{thremark}[prop]{Remark}
\newtheorem{thtable}[prop]{Table}

\newenvironment{defin}{\begin{thdefin}\rm}{\end{thdefin}}

\newenvironment{proof}{\medskip\noindent%
\textsc{Proof.}}{\qed\medskip}

\newenvironment{equationname}[1]{\settheequation{#1}\begin{equation}}%
{\end{equation}\restoretheequation}

\def\DJ{\leavevmode\setbox0=\hbox{D}\kern0pt\rlap
 {\kern.04em\raise.188\ht0\hbox{-}}D}

\begin{document}

\title{Some remarks on Kurepa's left factorial}
\author{Bernd C. Kellner}
\date{}
\maketitle

\abstract{We establish a connection between the subfactorial function $S(n)$ and the left factorial function
of Kurepa $K(n)$. Some elementary properties and congruences of both functions are described. Finally, we give a
calculated distribution of primes below 10000 of $K(n)$.}
\smallskip

\smallskip


\parindent 0cm

\section{Introduction}
\setcounter{equation}{0}

The subfactorial function is defined by
\[
   S(n) = n! \, \sum_{k=0}^n \frac{(-1)^k}{k!} \,, \quad n \in \NN_0
\]
which gives the number of permutations of $n$ elements without
any fixpoints, also called derangements of $n$ elements, see \cite[p.~195]{graham94concrete}.
This was already proven by P.~R.\ de~Montmort \cite{montmort08} in 1713.
L.\ Euler \cite{euler53} independently gave a proof in 1753, see also \cite{euler11}.
This function has the properties ($e \approx 2.71828$ is Euler's number)
\begin{eqnarray}
   S(n) &=& n S(n-1) + (-1)^n \,, \label{eqn-sn-1} \\
   S(n) &=& (n-1) \left( S(n-1) + S(n-2) \right) \,, \label{eqn-sn-2} \\
   S(n) &=& \left\lfloor \frac{n!}{e} \right\rfloor + \delta_n
      \quad \mbox{with} \quad \delta_n = \left\{
     \begin{array}{rl}
       0 \,, & 2 \notdiv n \\
       1 \,, & 2 \pdiv n
     \end{array} \right. \,. \label{eqn-sn-3}
\end{eqnarray}

Kurepa's left factorial function is defined by
\[
   K(0) = 0 \,, \qquad K(n) = \sum_{k=0}^{n-1} k! \,, \quad n \in \NN \,.
\]

In 1971 \DJ.\ Kurepa \cite{kurepa71} introduced the left factorial function
which is denoted by $!n = K(n)$.
Sometimes the subfactorial function is also denoted by $!n$, so we do not
use this notation to avoid confusion. For more details of the following conjecture
see a overview of A.\ Ivi\'c and \v{Z}.\ Mijajlovi\'c \cite{ivic95}.

\begin{conj}[Kurepa's left factorial hypothesis] \ \par
The following equivalent statements hold
\begin{equationname}{KH}
   \begin{array}{rcl@{\qquad}l}
   (K(n),n!) & = & 2 \,, & n \geq 2 \,, \\
   K(n) & \not\equiv & 0 \pmod{n} \,, & n > 2 \,, \\
   K(p) & \not\equiv & 0 \pmod{p} \,, & p \mbox{\ odd prime.}
   \end{array} \label{eqn-kh}
\end{equationname}
\end{conj}
\medskip

Recently, D.\ Barsky and B.\ Benzaghou \cite{barsky04} have given a proof of this
hypothesis. Since $K(n)$ is also related to Bell numbers $\Bell_n$ via
\[
   K(p) \equiv \Bell_{p-1} - 1 \pmod{p}
\]
for any prime $p$, they actually proved that $\Bell_{p-1} \not\equiv 1 \ \mod{p}$
is always valid for any odd prime $p$.
\medskip

Gessel \cite[Sect.~7/10]{gessel03} gives some recursive identities of
$S(n)$, $\Bell_n$, and others
with umbral calculus. Define symbolically $S^n = S(n)$ and $\Bell^n = \Bell_n$ with
$S^0 = \Bell^0 = 1$, then one may write
\begin{equation} \label{eqn-bs-umbral}
   \Bell^{n+1} = (\Bell+1)^n \quad \mbox{and} \quad
     n! = (S+1)^n \,, \quad n \geq 0 \,.
\end{equation}
Interestingly, both sequences have the same property as follows.
\begin{lemma}
Let $p$ be a prime. Then
\[
   \sum_{k=0}^p (-1)^k \, \Bell_k \equiv
     \sum_{k=0}^p (-1)^k \, S(k) \equiv 0 \pmod{p}
\]
with
\[
   \Bell_p \equiv 2 \pmod{p} \quad \mbox{and} \quad
     S(p) \equiv -1 \pmod{p} \,.
\]
\end{lemma}

\begin{proof}
By \refeqn{eqn-sn-1} and Wilson's theorem, we have $S(p) \equiv -1 \equiv (p-1)! \ \mod{p}$.
Hence, we can rewrite \refeqn{eqn-bs-umbral} by
$\Bell^{p} \equiv (\Bell+1)^{p-1}$ and $S^p \equiv (S+1)^{p-1} \ \mod{p}$.
Since $\binom{p-1}{k} \equiv (-1)^k \ \mod{p}$ for $0 \leq k < p$, this provides the
proposed congruence. Now, we use a congruence of Touchard for Bell numbers, see
\cite[Sect.~10, Theorem 10.1]{gessel03}. Then
\[
   \Bell_{n+p} - \Bell_{n+1} - \Bell_n \equiv 0 \pmod{p} \,, \quad n \geq 0 \,.
\]
With $n=0$ and $\Bell_0 = \Bell_1 = 1$, we obtain $\Bell_p \equiv 2 \ \mod{p}$.
\end{proof}
\medskip

First values of $K(n)$, $S(n)$, and $\Bell_n$ are given in the following table.
\newcommand{\mc}[1]{\multicolumn{1}{c|}{#1}}
\begin{center}
\begin{tabular}{|c*{12}{|r}|} \hline
  $n$    & \mc{0} & \mc{1} & \mc{2} & \mc{3} & \mc{4} & \mc{5} &
           \mc{6} & \mc{7} & \mc{8} & \mc{9} & \mc{10} \\ \hline \hline
  $K(n)$ & 0 & 1 & 2 & 4 & 10 & 34 & 154 &  874 &  5914 &  46234 & 409114 \\ \hline
  $S(n)$ & 1 & 0 & 1 & 2 &  9 & 44 & 265 & 1854 & 14833 & 133496 &1334961 \\ \hline
  $\Bell_n$ & 1 & 1 & 2 & 5 & 15 & 52 & 203 & 877 & 4140 & 21147 & 115975 \\ \hline
\end{tabular}
\end{center}

\section{Congruences between $K(n)$ and $S(n)$}
\setcounter{equation}{0}

\begin{lemma} \label{theor-kn-sn-1}
Let $n$ be a positive integer, then
\[
   K(n) \equiv (-1)^{n-1} S(n-1) \pmod{n} \,.
\]
\end{lemma}

\begin{proof}
Case $n=1$ is trivial. Let $n \geq 2$. Then we have
\[
   (-1)^{n-1} S(n-1) = \sum_{k=0}^{n-1} (-1)^{n-1-k} \binom{n-1}{k} (n-1-k)!
     = \sum_{k=0}^{n-1} (-1)^{k} \binom{n-1}{k} k!
\]
by turning the summation. Since it is valid for $0 \leq k < n$
\[
   (-1)^{k} \binom{n-1}{k} k! = (-1)^{k} (n-1) \cdots (n-k) \equiv k! \pmod{n} \,,
\]
this provides, term by term, the congruence claimed above.
\end{proof}

By Lemma \ref{theor-kn-sn-1} and \refeqn{eqn-sn-3},
we easily obtain a generalization, however, which is only noted for primes
elsewhere.

\begin{corl}
Let $n$ be a positive integer, then
\[
    K(n) \equiv (-1)^{n-1} \left\lfloor \frac{(n-1)!}{e} \right\rfloor + \delta_{n-1} \pmod{n} \,.
\]
\end{corl}

Hence, \refeqn{eqn-kh} is equivalent to
\[
    \left\lfloor \frac{(n-1)!}{e} \right\rfloor \not\equiv - \delta_{n-1} \pmod{n} \,,
      \quad n > 2 \,,
\]
while by recursive property \refeqn{eqn-sn-1}
\[
    \left\lfloor \frac{n!}{e} \right\rfloor \equiv - \delta_{n-1} \pmod{n} \,,
      \quad n \geq 1
\]
is always valid.

\begin{corl} Let $n$ be a positive integer, then
\refeqn{eqn-kh} is equivalent to
\[
   \left\lfloor \frac{n!}{e} \right\rfloor - \left\lfloor \frac{(n-1)!}{e} \right\rfloor
     \equiv 0 \pmod{n}  \quad \Longleftrightarrow \quad n = 1, 2 \,.
\]
\end{corl}

\begin{lemma} \label{lem-diff-Kpl}
Let $p$ be a prime, then
\[
   K(p) - K(p-l) \equiv - \frac{S(l-1)}{(l-1)!} \pmod{p} \,, \qquad l = 1,\ldots,p \,.
\]
\end{lemma}

\begin{proof}
Let $l \in \{ 1,\ldots,p \}$. We then have
\[
   K(p) - K(p-l) = \sum_{k=p-l}^{p-1} k! = \sum_{k=1}^{l} (p-k)!
     \equiv \sum_{k=1}^{l} \frac{(-1)^k}{(k-1)!} = - \frac{S(l-1)}{(l-1)!} \pmod{p} \,,
\]
since
\begin{equation} \label{eqn-diff-Kpl}
   (p-k)! \equiv \frac{(-1)^k}{(k-1)!} \pmod{p}
\end{equation}
follows by Wilson's theorem.
\end{proof}

\begin{corl} \label{corl-kh-ks}
Let $p$ be an odd prime, then \refeqn{eqn-kh} implies for $0 \leq l < p$
\[
   K(p-1-l) \not\equiv  \frac{S(l)}{l!} \pmod{p}
\]
respectively
\[
    l! \, K(p-1-l) \not\equiv \left\lfloor \frac{l!}{e} \right\rfloor + \delta_l \pmod{p} \,.
\]
\end{corl}
\medskip

Since \refeqn{eqn-kh} is true, we obtain, as an example, the following congruences
\[
   K(p) \not\equiv 0 \,, \  K(p-1) \not\equiv 1 \,, \  K(p-2) \not\equiv 0 \,, \
   K(p-3) \not\equiv \frac{1}{2} \,, \  K(p-4) \not\equiv \frac{1}{3}  \pmod{p} \,.
\]

\section{Properties of $K(n)$}
\setcounter{equation}{0}

To describe some interesting properties of $K(n)$,
we introduce the following definition which we name after Kurepa.

\begin{defin}
Let $p$ be an odd prime. The pair $(p,n)$ is called a \textsl{Kurepa pair}
if $p^r \pdiv K(n)$ with some integer $r \geq 1$. The max.\ integer $r$ is called
the \textsl{order} of $(p,n)$.
The \textsl{index} of $p$ is defined by
\[
   i_K(p) = \#\{ n \ : \ (p,n) \mbox{\ is a Kurepa pair} \} \,.
\]
If $i_K(p) > 0$, then $p$ is called a \textsl{Kurepa prime}.
\end{defin}
\medskip

We have, e.g., the Kurepa pairs $(19,7)$, $(19,12)$, and $(19,16)$.
If \refeqn{eqn-kh} would fail at an odd prime $p$, then this would imply $i_K(p) = \infty$.
This is an easy consequence of
\[
   p \pdiv K(p) \,, \quad p \pdiv (p+m)! \quad \mbox{for \ } m \geq 0 \,.
\]
The case $p=2$ is handled separately. One easily sees that $2 \pdiv K(n)$ for $n \geq 2$ and
$K(n) \equiv 2 \ \mod{4}$ for $n \geq 4$.
First values of $i_K(p)$ are given in the following table.
\medskip

\begin{center}
\begin{small}
\begin{tabular}{|c*{18}{|r}|} \hline
  $p$ & 3 & 5 & 7 & 11 & 13 & 17 & 19 & 23 & 29 & 31 & 37 & 41 & 43 & 47 & 53 & 59 & 61 \\ \hline \hline
  $\!i_K(p)\!$ & 0 & 1 & 1 & 1 & 0 & 1 & 3 & 1 & 0 & 2 & 1 & 2 & 0 & 0 & 0 & 0 & 1 \\ \hline
\end{tabular}
\end{small}
\end{center}
\smallskip

\begin{theorem}
Let $(p,n)$ be a Kurepa pair. Then
$p > n > 3$ is valid with
\[
   K(p) \equiv (-1)^n n! \, S(p-1-n) \pmod{p}
\]
which implies $p \notdiv S(p-1-n)$. Furthermore one has $i_K(p) \leq \lfloor \frac{p-1}{4} \rfloor$.
Consequently, there exist infinitely many Kurepa primes.
\end{theorem}

\begin{proof}
For now, let $p$ be an odd prime.
Let $(p,n)$ be a Kurepa pair. Since $p \notdiv K(p+m)$ for $m \geq 0$
by validity of \refeqn{eqn-kh}
and first values of $K(\cdot)$ are $0,1,2,4$, this yields $p > n > 3$.
We use Lemma \ref{lem-diff-Kpl} with $n=p-l$, then we have
\[
   0 \not\equiv K(p) \equiv K(p) - K(n) \equiv - \frac{S(p-1-n)}{(p-1-n)!} \pmod{p}
\]
which provides the result by means of \refeqn{eqn-diff-Kpl} and also
$p \notdiv S(p-1-n)$. Now, we have to count possible Kurepa pairs.
Corollary \ref{corl-kh-ks} shows that $K(p-2) \not\equiv 0 \ \mod{p}$.
If $p \pdiv K(n)$ then $p \notdiv K(n+l)$ for $l=1,2,3$. This is seen by
$n! \not\equiv 0 \ \mod{p}$ and
\[
   n!+(n+1)! = (n+2) \, n! \not\equiv 0 \,, \
     n!+(n+1)!+(n+2)! = (n+2)^2 \, n! \not\equiv 0 \pmod{p} \,,
\]
since $n \neq p-2$. On the other side, we have $4 \leq n \leq p-1$. Then a simple
counting argument provides $i_K(p) \leq \lfloor \frac{p-1}{4} \rfloor$.
Finally, $K(n) \rightarrow \infty$ for $n \rightarrow \infty$ and
$p \pdiv K(n) \Rightarrow p>n$ for odd primes
imply infinitely many Kurepa primes.
\end{proof}
\medskip

Now, the remarkable fact of $K(n)$ is the finiteness of Kurepa pairs for all odd primes.
In $p$-adic analysis, the series
\[
    K(\infty) = \sum_{k=0}^\infty k!
\]
is an example of a convergent series resp.\ $K(n)$ is a convergent sequence which
lies in $\ZZ_p$, the ring of $p$-adic integers. Then \refeqn{eqn-kh} is equivalent
to $K(\infty)$ is a unit in $\ZZ_p$ for all odd primes $p$. The behavior $\mod{p^r}$
is illustrated by the following theorem. Note that
$l_r$ is related to the so-called Smarandache function for factorials.

\begin{theorem}
Let $p, r$ be positive integers with $p$ prime. Then the sequence
\[
    K(n) \ \mod{p^r} \,, \quad n \geq 0
\]
is constant for $n \geq l_r p$ with $r \geq l_r$ and
\[
   l_r = \min_l \left\{ l \ : \ l + \frac{l-\sigma_p(l)}{p-1} \geq r \right\} \,,
\]
where $\sigma_p(l)$ gives the sum of digits of $l$ in base $p$.
\end{theorem}

\begin{proof} We have to determine a minimal $l$ with the property $\ord_p \, (lp)! \geq r$.
Counting factors which are divisible by $p$, we obtain
\[
    \ord_p \, (lp)! = l + \ord_p l! = l + \frac{l-\sigma_p(l)}{p-1}
\]
by means of the $p$-adic valuation of factorials, see \cite[Section 3.1, p.~241]{robert00padic}.
\end{proof}
\pagebreak

At the end, we give some results of calculated Kurepa pairs.
There are $N = \pi(10000)-1 = 1228$ odd primes below 10000.
Let $N_r$ be the number of odd primes with index $i_K(p)=r$ in this range.
The following table shows the distribution of the index $i_K$.

\begin{center}
\begin{tabular}{|c*{7}{|c}|} \hline
 $r$     & 0 & 1 & 2 & 3 & 4 & 5 \\ \hline \hline
 $N_r$   & 459 & 472 & 213 & 58 & 23 & 3 \\ \hline
 $N_r/N$ & 0.37378 & 0.38436 & 0.17345 & 0.04723 & 0.01873 & 0.00244 \\ \hline
\end{tabular}
\end{center}
\medskip

The calculated Kurepa pairs with index $i_K(p)=5$ are as follows.
\begin{center}
\begin{small}
\begin{tabular}{*{5}{|c}|} \hline
(2203,277)  & (2203,788)  & (2203,837)  & (2203,1246) & (2203,1927) \\ \hline
(5227,850)  & (5227,1752) & (5227,3451) & (5227,4363) & (5227,4716) \\ \hline
(6689,1716) & (6689,2404) & (6689,3641) & (6689,3969) & (6689,6601) \\ \hline
\end{tabular}
\end{small}
\end{center}
\medskip

All primes below 10000 appear with a simple power in $K(n)$, except $K(3)=4$.
On the other side, the occurrence of higher powers $p^r$ in $K(n)$ seems to be very rare.
M.~Zivkovic \cite{ziv99} gives the first example $54503^2 \pdiv K(26541)$.
There are two Kurepa pairs $(54503,26541)$ and $(54503,49783)$, but only
the first of them has order two.
\medskip

One may ask whether the distribution of Kurepa pairs resp.\ the index $i_K$ can be
asymptotically determined and even proven. Are there infinitely many non-Kurepa primes
$p$ with $i_K(p)=0$? It seems that this subject of $K(n)$ and its distribution
of primes will be much simpler to attack as, for example,
the more complicated but in a sense similar case of
the distribution of irregular primes of Bernoulli numbers.

\subsubsection*{Acknowledgement}
The author wishes to thank Prof.\ Ivi\'c for informing about the problem of Kurepa.
\medskip

Bernd C. Kellner \\
{\small
address: Reitstallstr. 7, 37073 G\"ottingen, Germany \\
email: bk@bernoulli.org}

\bibliographystyle{plain}
\bibliography{ku_bib}

\end{document}